\documentclass[12pt]{amsart}

\usepackage{mathrsfs}

\usepackage{amsmath,amsthm,amssymb}
\font\teneufm=eufm10
\font\seveneufm=eufm7
\font\fiveeufm=eufm5
\newfam\frakturfam

\textfont\frakturfam=\teneufm
\scriptfont\frakturfam=\seveneufm
\scriptscriptfont\frakturfam=\fiveeufm

\newtheorem{pr}{Proposition}

\newtheorem{lm}{Lemma}
\newtheorem{theor}{Theorem}
\newtheorem{co}{Corollary}

\newtheorem{prob}{Problem}
\def\bee{\begin{eqnarray}}
\def\bes{\begin{eqnarray*}}
\def\eee{\end{eqnarray}}
\def\ees{\end{eqnarray*}}

\def\Proof{{\sl Proof.}\ }

\title{Associative, Lie, and left-symmetric algebras of derivations}

\begin{document}
\date{}
\maketitle

\begin{center}

{\bf Ualbai Umirbaev}\footnote{Supported by an NSF grant DMS-0904713 and by an MES grant 0755/GF of Kazakhstan; Eurasian National University,
 Astana, Kazakhstan and
 Wayne State University,
Detroit, MI 48202, USA,
e-mail: {\em umirbaev@math.wayne.edu}}

\end{center}

\begin{abstract} Let $P_n=k[x_1,x_2,\ldots,x_n]$ be the polynomial algebra over a field $k$ of characteristic zero in the variables $x_1,x_2,\ldots,x_n$ and $\mathscr{L}_n$ be the left-symmetric algebra of all derivations of $P_n$ \cite{Dzhuma99,UU2014-1}.  Using the language of $\mathscr{L}_n$, for every derivation $D\in \mathscr{L}_n$ we define the associative algebra $A_D$, the Lie algebra $L_D$, and the left-symmetric algebra $\mathscr{L}_D$ related to the study of the Jacobian Conjecture. For every derivation $D\in \mathscr{L}_n$ there is a unique $n$-tuple $F=(f_1,f_2,\ldots,f_n)$ of elements of $P_n$ such that $D=D_F=f_1\partial_1+f_2\partial_2+\ldots+f_n\partial_n$. In this case, using an action of the Hopf algebra of noncommutative symmetric functions $\mathrm{NSymm}$ on $P_n$, we show that these algebras are closely related to the description of coefficients  of the formal inverse to  the polynomial endomorphism $X+tF$, where $X=(x_1,x_2,\ldots,x_n)$ and $t$ is an independent parameter.

We prove that the Jacobian matrix $J(F)$ is nilpotent if and only if all right powers $D_F^{[r]}$ of $D_F$ in $\mathscr{L}_n$ have zero divergence. In particular, if $J(F)$ is nilpotent then $D_F$ is right nilpotent.

We discuss some advantages and shortcomings of these algebras and formulate some open questions.
\end{abstract}


\noindent


\section{Introduction}

\hspace*{\parindent}

Let $k$ be an arbitrary field of characteristic zero  and $P_n=k[x_1,x_2,\ldots,x_n]$ be the polynomial algebra over $k$ in the variables $x_1,x_2,\ldots,x_n$. There are two well known algebras related to the study of derivations of $P_n$. They are the Witt algebra $W_n$ and the Weyl algebra $A_n$. Recall that $W_n$ is the Lie algebra of all derivations of $P_n$ and $A_n$ is the associative algebra of all linear differential operators on $P_n$.

The set of elements $u\partial_i$, where $u=x_1^{s_1}\ldots x_n^{s_n}$ is an arbitrary monomial, $\partial_i=\frac{\partial}{\partial x_i}$, and $1\leq i\leq n$, forms a linear basis for $W_n$. For any $u=a\partial_i, v=b\partial_j$, where $a,b\in P_n$ are monomials, put
\bee\label{Yf1}
u\cdot v= ((a\partial_i)(b))\partial_j.
\eee
Extending this operation by distributivity, we get a well defined bilinear operation $\cdot$ on $W_n$.
Denote this algebra by $\mathscr{L}_n$. It is easy to check (see Section 2) that $\mathscr{L}_n$ is a left-symmetric algebra \cite{Dzhuma99,UU2014-1} and its commutator algebra is the Witt algebra $W_n$. We say that $\mathscr{L}_n$ is {\em the left-symmetric algebra of derivations} of $P_n$.

The language of the left-symmetric algebras of derivations is very convenient to describe some important notions of affine algebraic geometry in purely algebraic terms \cite{UU2014-1}. For example, an element of $\mathscr{L}_n$ is left nilpotent if and only if it is a locally nilpotent derivation of $P_n$. One of the greatest algebraic advantages of $\mathscr{L}_n$ is that $\mathscr{L}_n$ satisfies an exact analogue of the Cayley-Hamilton trace identity. Recall that $W_n$ and $A_n$ do not have an analogue of this identity.

Let $D\in \mathscr{L}_n$ be an arbitrary derivation of $P_n$. Denote by $\mathscr{L}_D$ the subalgebra of the left-symmetric algebra $\mathscr{L}_n$ generated by $D$. Denote by $L_D$ the Lie subalgebra of the Witt algebra $W_n$ generated by all right powers $D^{[p]}$ of $D$. Obviously, $L_D\subseteq \mathscr{L}_D$. Denote by $A_D$ the subalgebra (with identity) of the Weyl algebra $A_n$ generated by all right powers $D^{[p]}$ of $D$. So, $A_D$ is an associative enveloping algebra of the Lie algebra $L_D$. The Lie algebra $L_D$ is a nontrivial Lie algebra ever related to one derivation.

Every $n$-tuple $F=(f_1,f_2,\ldots,f_n)$ of elements of $P_n$ represents a polynomial endomorphism of the vector space $k^n$. We denote by $F^*$ the endomorphism of $P_n$ defined by $F^*(x_i)=f_i$ for all $i$. Also denote  by
\bes
D_F=f_1\partial_1+f_2\partial_2+\ldots+f_n\partial_n
\ees
the derivation of $P_n$ defined by $D_F(x_i)=f_i$ for all $i$. Note that every derivation $D$ can be uniquely represented as $D=D_F$ for some polynomial $n$-tuple $F$.  Using this correspondence we often use parallel notations  $A_D=A_F$, $L_D=L_F$, and $\mathscr{L}_D=\mathscr{L}_F$ if $D=D_F$

We show that if the Jacobian matrix $J(F)$ is nilpotent then $D_F$ is a right nilpotent element of $\mathscr{L}_n$.
We also show that the Jacobian matrix $J(F)$ is nilpotent if and only if all  right powers $D_F^{[p]}$ of $D_F$ have zero divergence. Moreover, if $J(F)$ is nilpotent then every element of $L_F$ has zero divergence.

Let $t$ be an independent parameter and
\bes
(X+tF)^{-1}=X+tF_1+t^2F_2+\ldots+t^nF_n+\ldots
\ees
be the formal (or analytic) inverse to the endomorphism $X+tF$ of $k[t]^n$. There are many interesting papers devoted to the description of $F_i$ \cite{BCW,Essen,Wright07}. We show that $A_F$ and $L_F$ are also generated by all $D_{F_i}$ where $i\geq 1$. For this reason we can say that $A_F$ and $L_F$ are, respectively, the associative and the Lie algebras of coefficients of the formal inverse to $X+tF$. Notice that $\mathscr{L}_F$ is also the smallest left-symmetric algebra containing all $D_{F_i}$ where $i\geq 1$.

Recall that the Hopf algebra of noncommutative symmetric functions $\mathrm{NSymm}$ \cite{Gelfand} regarded as an algebra  is the free associative algebra
\bes
\mathrm{NSymm}=k\langle Z_1,Z_2,\ldots,Z_n,\ldots\rangle
\ees
over $k$ in the variables $Z_1,Z_2,\ldots,Z_n,\ldots$.

We define an action of $\mathrm{NSymm}$ on $P_n$ by
\bes
(X+tF)^*(a)=a+tZ_1(a)+t^2Z_2(a)+\ldots+t^nZ_n(a)+\ldots
\ees
for any $a\in P_n$. This action represents a natural linearization of the action of $(X+tF)^*$ on $P_n$.  We show that $A_F$ is the image of $\mathrm{NSymm}$ under this representation and $L_F$ is the image of the Lie algebra $\mathrm{Prim}$ of all primitive elements of $\mathrm{NSymm}$. In this way, $A_F$ and $L_F$ may be considered as linearization algebras of the action of $(X+tF)^*$ on $P_n$. The left-symmetric algebra $\mathscr{L}_F$ also can be related to further linearizations.

The Hopf algebra of noncommutative symmetric functions $\mathrm{NSymm}$ was introduced in \cite{Gelfand} as a noncommutative generalization of the Hopf algebra of symmetric functions $\mathrm{Symm}$. Several systems of free and primitive generators of $\mathrm{NSymm}$ and relations between them were given in \cite{Gelfand}.  Some more relations between the generators of $\mathrm{NSymm}$ are given in \cite{Zhao08}.

There are two well known systems of free primitive generators  \cite{Hazewinkel05} of $\mathrm{NSymm}$ which are dual to each other with respect to the standard involution of the free associative algebra $k\langle Z_1,Z_2,\ldots,Z_n,\ldots\rangle$.  It is interesting that one of them corresponds to the right powers $D_F^{[r]}$ of $D_F$ and the other one corresponds to the $D_{F_i}$ for all $i\geq 1$. These observations make the Lie algebra $L_F$ very important in studying the Jacobian Conjecture. The right powers $D_F^{[r]}$ are very convenient to express that $J(F)$ is nilpotent. In order to solve the Jacobian Conjecture it is necessary to prove that there exists a positive integer $m$ such that $D_{F_i}=0$ for all $i\geq m$.

The action of $\mathrm{NSymm}$ on $P_n$, defined above, corresponds to one of a series of homomorphisms constructed in  \cite{Zhao08} and the images of primitive generators were calculated in \cite{Zhao08}.

It is rewarding to initiate a systematic study of the associative algebra $A_F$, the Lie algebra $L_F$ and the left-symmetric algebra $\mathscr{L}_F$.
Using an example of an automorphism studied earlier by A. van den Essen \cite{Essen85} and  G. Gorni and G. Zampieri \cite{GZ96}, we give an example of $F$ with nilpotent Jacobian matrix $J(F)$ such that $L_F$ is not nilpotent nor solvable.

The paper is organized as follows. Section 2 is devoted to the study of the left-symmetric algebra $\mathscr{L}_n$. In particular, we describe the right and the left multiplication algebras of $\mathscr{L}_n$ and describe an analogue of the Cayley-Hamilton identity. In Section 3 we develop technics for calculation of divergence of elements in $\mathscr{L}_n$. The definition of the Hopf algebra of noncommutative symmetric functions $\mathrm{NSymm}$ is given in Section 4. We give also some primitive systems of generators of $\mathrm{NSymm}$ and relations between from \cite{Gelfand}. The action of $\mathrm{NSymm}$ and the images of primitive elements are given in Section 5. In Section 6 we discuss some properties of these algebras towards the Jacobian Conjecture and formulate some open problems.

\section{Algebra $\mathscr{L}_n$}

\hspace*{\parindent}

If $A$ is an arbitrary linear algebra over a field $k$ then the set $\mathrm{Der}_kA$ of all $k$-linear derivations  of $A$ forms a Lie algebra. If $A$ is a free algebra then it is possible to define a multiplication  $\cdot$ on $\mathrm{Der}_kA$ such that it becomes a left-symmetric algebra and its commutator algebra becomes the Lie algebra of derivations $\mathrm{Der}_kA$ of $A$ \cite{UU2014-1}.

Recall that an algebra $\mathscr{L}$  over $k$ is called left-symmetric \cite{Burde2006} if $\mathscr{L}$ satisfies the identity
\bee\label{Yf2}
(xy)z-x(yz)=(yx)z-y(xz).
\eee
This means that the associator $(x,y,z):=(xy)z-x(yz)$ is symmetric with respect to two left arguments, i.e.,
\bes
(x,y,z)=(y,x,z).
\ees
The variety of
left-symmetric algebras is Lie-admissible, i.e., each
left-symmetric algebra $\mathscr{L}$ with the operation $[x,y]:=xy-yx$ is a
Lie algebra.

Recall that the space of the algebra $\mathscr{L}_n$ is $W_n$ and the product is defined by (\ref{Yf1}).
\begin{lm}\label{Yl1} \cite{Dzhuma99,UU2014-1} Algebra $\mathscr{L}_n$ is left-symmetric and its commutator algebra is the Witt algebra $W_n$.
\end{lm}
\Proof Let $x,y\in \mathscr{L}_n$. Denote by $[x,y]=x\cdot y-y\cdot x$ the commutator of $x$ and $y$ in $\mathscr{L}_n$ and denote by $\{x,y\}$ the product of $x$ and $y$ in $W_n$. We first prove that the commutator algebra of $\mathscr{L}_n$ is $W_n$, i.e.,
\bes
[x,y](a)=\{x,y\}(a)
\ees
 for all $a\in P_n$. Note that
\bes
\{x,y\}(a)=x(y(a))-y(x(a))
\ees
by the definition. Taking into account that $[x,y]$ and $\{x,y\}$ are both derivations, we can assume that $a=x_t$.
Consequently, it is sufficient to check that
\bes
(x\cdot y-y\cdot x)(x_t)=x(y(x_t))-y(x(x_t)).
\ees
We may also assume that $x=u\partial_i$ and $y=v\partial_j$. If $t\neq i,j$, then all components of the last equality are zeroes. If $t=i\neq j$ or $t=i=j$, then it is also true. Consequently, the commutator algebra of $\mathscr{L}_n$ is $W_n$.

Assume that $x,y\in \mathscr{L}_n$ and $z=a\partial_t$. Then
\bes
(x,y,z)=(xy)z-x(yz)=[(xy)(a)-x(y(a))]\partial_t,\\
(y,x,z)=(yx)z-y(xz)=[(yx)(a)-y(x(a))]\partial_t.
\ees
To prove (\ref{Yf2}) it is sufficient to check that
\bes
[x,y](a)=x(y(a))-y(x(a))=\{x,y\}(a),
\ees
which is already proved. $\Box$

A natural $P_n$-module structure on $\mathscr{L}_n$ can be defined by $p \cdot  u\partial_i=(pu)\partial_i$ for all $i$ and $p,u\in P_n$. Then
\bes
\mathscr{L}_n=P_n \partial_1\oplus P_n \partial_2\oplus\ldots \oplus P_n \partial_n
\ees
is a free $P_n$-module.

Consider the grading
\bes
P_n=A_0\oplus A_1\oplus A_2\oplus\ldots\oplus A_s\oplus\ldots,
\ees
where $A_i$ the space of homogeneous elements of degree $i\geq 0$. The left-symmetric algebra $\mathscr{L}_n$ has a natural grading
\bes
\mathscr{L}_n=L_{-1}\oplus L_{0}\oplus L_{1}\oplus \ldots \oplus L_{s}\oplus\ldots,
\ees
where $L_{i}$ the space of elements of the form $a\partial_j$ with $a\in A_{i+1}$ and $1\leq j\leq n$. Elements of $L_{s}$ are called homogeneous derivations of $P_n$ of degree $s$.

We have $L_{-1}=k\partial_1+\ldots+k\partial_n$ and $L_{0}$ is a subalgebra of $\mathscr{L}_n$ isomorphic to the matrix algebra $M_n(k)$. The element
\bes
D_X=x_1\partial_1+x_2\partial_2+\ldots+x_n\partial_n
\ees
 is the identity element of the matrix algebra $L_{0}$ and is the right identity element of $\mathscr{L}_n$. The left-symmetric algebra
$\mathscr{L}_n$ has no identity element.

We establish some properties of $\mathscr{L}_n$ related to the Jacobian Conjecture.

For every $n$-tuple $F=(f_1,f_2,\ldots,f_n)$ of elements of $P_n$ denote by $J(F)=(\partial_j(f_i))_{1\leq i,j\leq n}$  the Jacobian matrix of $F$. Notice that every derivation $D$ of $P_n$ has the form $D=D_F$ for a unique endomorphism $F$. Put $J(D)=J(F)$. So, the Jacobian matrix of every derivation $D$ of $A$ is defined.

\begin{lm}\label{Yl2} \cite{UU2014-1} Let $F$ and $G$ be two arbitrary $n$-tuples of elements of $A$.
Then
\bes
D_F D_G=D_{D_F(G)}=D_{J(G)F}=D_{J(D_G)F}.
\ees
\end{lm}
\Proof The definition of the left symmetric product $\cdot$ directly implies that $D_F D_G=D_{D_F(G)}$.
Notice that for any $h\in A$ we have
\bes
D_F(h)=\sum_{i=1}^n \frac{\partial h}{\partial x_i}y_i|_{y_i:=f_i}
=(\frac{\partial h}{\partial x_1},\ldots,\frac{\partial h}{\partial x_n})F.
\ees
Consequently, $D_{D_F(G)}=D_{J(G)F}$. $\Box$

 For any $a\in \mathscr{L}_n$ put $a^0=a^{[0]}=a$, $a^{r+1}=a(a^r)$, and $a^{[r+1]}=(a^{[r]})a$ for any $r\geq 0$. It is natural to say that $a$ is left nilpotent if $a^m=0$ for some $m\geq 2$. Similarly,  $a$ is right nilpotent if $a^{[m]}=0$ for some $m\geq 2$.

\begin{lm}\label{Yl3} \cite{UU2014-1}
A derivation $D$ of $A$ is locally nilpotent if and only if $D$ is a left nilpotent element of $\mathscr{L}_n$.
\end{lm}
\Proof Suppose that $D=D_F$ and put
\bes
H_i=\underbrace{D(D\ldots(D(D}_iX))\ldots)
\ees
 for all $i\geq 1$. Note that $H_1=F$ and $H_2=D_F(F)$.  Consequently, $D^2=D_{H_2}$  by Lemma \ref{Yl2}.
 Continuing the same calculations,  it is easy to show that $D^i=D_{H_i}$ for all $i$. Consequently, $D^m=0$ if and only if $H_m=0$. Note that $H_m=0$ means that $D$ applied $m$ times to $x_i$  gives $0$ for all $i$. $\Box$

{\bf Example 1.}
Consider a well known \cite{Bass2} locally nilpotent derivation
\bes
D=(x^2-yz)(z\frac{\partial}{\partial x}+2x\frac{\partial}{\partial y})
\ees
 of $k[x,y,z]$. It is easy to check that $D$ is not right nilpotent. So, the left nilpotency of derivations does not imply their right nilpotency.

Let $\mathscr{L}$ be an arbitrary left-symmetric algebra. Denote by $\mathrm{Hom}_k(\mathscr{L},\mathscr{L})$ the associative algebra of all $k$-linear transformations of the vector space $\mathscr{L}$.  For any $x\in \mathscr{L}$ denote by $L_x : \mathscr{L}\rightarrow \mathscr{L} (a\mapsto xa)$ and $R_x : \mathscr{L}\rightarrow \mathscr{L} (a\mapsto ax)$ the operators of left and right multiplication by $x$, respectively. It follows from (\ref{Yf2}) that
\bee\label{Yf3}
L_{[x,y]}=[L_x,L_y],\ \ \
R_{xy}=R_yR_x+[L_x,R_y].
\eee
Denote by $M(\mathscr{L})$ the subalgebra of $\mathrm{Hom}_k(\mathscr{L},\mathscr{L})$ (with identity) generated by all $R_x,L_x$, where $x\in \mathscr{L}$. Algebra $M(\mathscr{L})$ is called the {\em multiplication} algebra of $\mathscr{L}$. The subalgebra $R(\mathscr{L})$ of $M(\mathscr{L})$ (with identity) generated by all $R_x$, where $x\in \mathscr{L}$, is called the {\em right multiplication} algebra of $\mathscr{L}$. Similarly, the subalgebra $L(\mathscr{L})$ of $M(\mathscr{L})$ (with identity) generated by all $L_x$, where $x\in \mathscr{L}$, is called the {\em left multiplication} algebra of $\mathscr{L}$.

\begin{lm}\label{Yl4} The right multiplication algebra $R(\mathscr{L}_n)$ of $\mathscr{L}_n$ is isomorphic to the matrix algebra $M_n(P_n)$ and there exists a unique isomorphism $\theta : R(\mathscr{L}_n)\rightarrow M_n(P_n)$ such that $\theta(R_D)=J(D)$ for all $D\in \mathscr{L}_n$.
\end{lm}
\Proof Let $D\in \mathscr{L}$. Notice that $R_D=0$ if and only $D\in  k\partial_1+\ldots+k\partial_n=L_0$. In fact, suppose that $R_D=0$. Then $\partial_i\cdot D=0$ for all $i$. This means that if $D=D_F$ then $F$ does not contain $x_i$ for all $i$ and $D\in L_0$. Consequently, $R_D=0$ if and only if $J(D)=0$.

Thus the correspondence $R_D\mapsto J(D)$ is well defined.
 Notice that for any $D=D_F,D_1,\ldots,D_m$ we have
\bes
R_{D_1}\ldots R_{D_m}(D)=(\ldots(D\cdot D_m)\ldots D_1)=D_{J(D_1)\ldots J(D_m)F}
\ees
by Lemma \ref{Yl2}. This implies that the equality $f(R_{D_1},\ldots,R_{D_m})=0$, where $f$ is an associative polynomial, holds if and only if
$f(J(D_1),\ldots,J(D_m))=0$. Consequently, there exists a unique monomorphism   $\theta : R(\mathscr{L})\rightarrow M_n(P_n)$ such that $\theta(R_D)=J(D)$ for all $D\in \mathscr{L}_n$. The uniqueness of $\theta$ is obvious since $R(\mathscr{L}_n)$ is generated by all $R_D$.

Denote by $B$ the subalgebra of $M_n(A)$ generated by all Jacobian matrices. Denote by $e_{ij}$, where $1\leq i,j\leq n$, the matrix with $1$ in the $(i,j)$ place and with zeroes everywhere else, i.e., the matrix identities. Consider $F=(f_1,\ldots,f_n)$.
 If $f_i=x_j$ and $f_s=0$ for all $s\neq i$ then $J(F)=e_{ij}$ and $e_{ij}\in B$ for all $i,j$.
 Let $u=x_1^{s_1}\ldots x_n^{s_n}$ be an arbitrary monomial of $P_n$. Put $f_1=1/(s_1+1) x_1^{s_1}x_2^{s_2}\ldots x_n^{s_n}$ and $f_i=0$ for all $i\geq 2$. Then $u$ becomes the element of $J(F)$ in the place $(1,1)$. This implies that $e_{i1}J(F)e_{1j}=ue_{ij}$. Consequently, $B=M_n(P_n)$ and $\theta$ is a surjection. $\Box$

Identities of $\mathscr{L}_n$ are studied by A.S. Dzhumadildaev \cite{Dzhuma99,Dzhuma00,Dzhuma04}. If $n=1$ then $\mathscr{L}_1$ becomes a Novikov algebra and identities of  $\mathscr{L}_1$ are studied in \cite{MLU2011-1}.

\begin{co}\label{Yc1}
The identities of the right multiplication algebra $R(\mathscr{L}_n)$ coinside with the identities of the matrix algebra $M_n(k)$.
\end{co}

\begin{co}\label{Yc2} \cite{UU2014-1}
Let $D\in\mathscr{L}_n$. Then the Jacobian matrix $J(D)$ of $D$ is nilpotent if and only if  $R_D$ is a nilpotent element of $M(\mathscr{L}_n)$.
\end{co}
\Proof By Lemma \ref{Yl4}, $J(D)^s=0$ if and only if $R_D^s=0$. $\Box$

Consequently, if $J(D)$ is nilpotent then $D$ is right nilpotent.
Is the converse true? This question is still open.

Every element $p\in P_n$ can be considered as an element of $\mathrm{Hom}(\mathscr{L}_n,\mathscr{L}_n)$ since $\mathscr{L}_n$. Then $P_n R(\mathscr{L}_n)$ becomes a left $P_n$-module. Notice that  $M_n(P_n)$ is also a $P_n$-module.

\begin{lm}\label{Yl5} $P_n R(\mathscr{L}_n)=R(\mathscr{L}_n)$ and the isomorphism $\theta : R(\mathscr{L}_n)\rightarrow M_n(P_n)$, constructed in Lemma \ref{Yl4}, is an isomorphism of $P_n$-modules.
\end{lm}
\Proof As in the proof of Lemma \ref{Yl4}, for any $p\in P_n$ and $D=D_F,D_1,\ldots,D_m$ we have
\bes
p R_{D_1}\ldots R_{D_m}(D)=(\ldots(D\cdot D_m)\ldots D_1)=D_{pJ(D_1)\ldots J(D_m)F}
\ees
by Lemma \ref{Yl2}. This implies that the equality $f(R_{D_1},\ldots,R_{D_m})=0$, where $f$ is an associative polynomial over $P_n$, holds if and only if
$f(J(D_1),\ldots,J(D_m))=0$. Consequently, there exists a unique monomorphism   $\overline{\theta} : P_n R(\mathscr{L}_n)\rightarrow M_n(P_n)$ of $P_n$-modules such that $\overline{\theta}(T)=\theta(T)$ for all $T\in R(\mathscr{L}_n)$. Then $\overline{\theta}$ is an isomorphism since $\theta$ is an isomorphism. This implies that
$P_n R(\mathscr{L}_n)=R(\mathscr{L}_n)$ and $\overline{\theta}=\theta$. $\Box$

The isomorphism $\theta : R(\mathscr{L}_n)\rightarrow M_n(P_n)$ from Lemma \ref{Yl4} gives us the matrix $\theta(T)$ for any $T\in R(\mathscr{L})$. Notice that $R_{D_X}$ is the identity element of $R(\mathscr{L}_n)$ and will be denoted  by $E$.

Let $T$  be an arbitrary element of $R(\mathscr{L})$. Then the matrix $\Theta=\theta(T)$ satisfies the well-known Cayley-Hamilton identity
\bes
\Theta^n+a_1\Theta^{n-1}+\ldots+a_{n-1}\Theta+a_nI=0,
\ees
 where $I$ is the identity matrix of order $n$ and $a_i\in P_n$. Recall that $a_1,a_2,\ldots,a_n$ can be expressed by traces of powers of $J$. It follows that
\bee\label{Yf4}
T^n+a_1T^{n-1}+\ldots+a_{n-1}T+a_nE=0
\eee
since $\theta$ is an isomorphism. This identity is an analogue of the Cayley-Hamilton trace identity for $\mathscr{L}$. Notice that if $T=f(R_{D_1},\ldots,R_{D_m})$ then $\Theta=f(J(D_1),\ldots,J(D_m))$. So, all coefficients of (\ref{Yf4}) can be expressed by traces of products of Jacobian matrices.

Yu. Razmyslov proved \cite{Razmyslov74} that all trace identities (in particular, all identities) of the matrix algebra $M_n(k)$ are corollaries of the Cayley-Hamilton trace identity. Consequently, all identities of $R(\mathscr{L}_n)$ are corollaries of (\ref{Yf4}).
Of course, every identity of $R(\mathscr{L}_n)$ gives a right identity of $\mathscr{L}_n$, i.e., an identity of
$\mathscr{L}_n$ which can be expressed by right multiplication operators. But it does not mean that every right multiplication operator identity of $\mathscr{L}_n$ is an identity of $R(\mathscr{L}_n)$. For this reason, we cannot say that every right identity of $\mathscr{L}_n$ is a corollary of (\ref{Yf4}).

\begin{lm}\label{Yl6} The left multiplication algebra  $L(\mathscr{L}_n)$ of $\mathscr{L}_n$ is isomorphic to the Weyl algebra $A_n$.
\end{lm}
\Proof Notice that for any $D=D_F,D_1,D_2,\ldots,D_m$ we have
\bes
L_{D_1}L_{D_2}\ldots L_{D_m}(D)=(D_1\ldots(D_m\cdot D)\ldots)=D_{D_1(D_2(\ldots D_m(F)\ldots))}.
\ees
This implies that the equality $f(L_{D_1},L_{D_2},\ldots,L_{D_m})=0$, where $f$ is an associative polynomial, holds in $L(\mathscr{L}_n)$ if and only if
$f(D_1,D_2,\ldots,D_m)=0$ holds in $A_n$. Consequently, there exists a unique monomorphism   $\psi : L(\mathscr{L}_n)\rightarrow A_n$ such that $\psi(L_D)=D$ for all $D\in \mathscr{L}_n$. Then $\psi$ is an epimorphism since $A_n$ is generated by all derivations. $\Box$

So, Lemmas \ref{Yl4} and \ref{Yl6} describe the structure of the right and left multiplicative algebras of $\mathscr{L}_n$, respectively. But at the moment I do not know the structure of the multiplication algebra $M(\mathscr{L}_n)$.
Recall that the Weyl algebra $A_n$ does not satisfy any nontrivial identity.  The left operator identities of $\mathscr{L}_n$ are very important in studying the locally nilpotent derivations and the Jacobian Conjecture.

\begin{lm}\label{Yl7} Let $f=f(z_1,z_2,\ldots,z_t)$ be a Lie polynomial. Then $f(z_1,z_2,\ldots,z_t)=0$  is an identity of the Witt algebra $W_n$ if and only if $f(L_{z_1},L_{z_2},\ldots,L_{z_t})=0$ is a left operator identity
of $\mathscr{L}_n$.
\end{lm}
\Proof Let $w_1,w_2,\ldots,w_t\in W_n=\mathscr{L}_n$. Notice that $f(w_1,w_2,\ldots,w_t)=0$ in $W_n$ if and only if $L_{f(w_1,w_2,\ldots,w_t)}=0$ in $L(\mathscr{L}_n)$ since the left annihilator of $\mathscr{L}_n$ is trivial. By (\ref{Yf3}), we get
\bes
L_{f(w_1,w_2,\ldots,w_t)}=f(L_{w_1},L_{w_2},\ldots,L_{w_t})=0.
\ees
This means that the associative polynomial
\bes
L_f=f(L_{z_1},L_{z_2},\ldots,L_{z_t})
\ees
in $L_{z_1},L_{z_2},\ldots,L_{z_t}$ is a  left operator identity of $\mathscr{L}_n$ if and only if $f(z_1,z_2,\ldots,z_t)$ is an identity of $W_n$.

Identities of $W_n$ are studied in \cite{Razmyslov94} and left operator identities of $\mathscr{L}_n$ are studied in \cite{Dzhuma04}.

\section{Divergence calculations}

\hspace*{\parindent}

If $D$ is an arbitrary element of $\mathscr{L}_n$, then there exists a unique $n$-tuple $F=(f_1,f_2,\ldots,f_n)$ of elements of $P_n$ such that  $D=D_F\in\mathscr{L}_n$. Put
\bes
\mathrm{div}(D)=\mathrm{div}(D_F)=\partial_1(f_1)+\partial_2(f_2)+\ldots+\partial_n(f_n).
\ees
Consequently, $\mathrm{div}(D)=\mathrm{Tr}(J(D))=\mathrm{Tr}(J(F))$.

Recall that every $n$-tuple $F=(f_1,f_2,\ldots,f_n)$ of $P_n$ represents a polynomial mapping of the vector space $k^n$. Denote by $F^*$  the endomorphism of $P_n$ such that $F^*(x_i)=f_i$ for all $i$.
If $F$ and $G$ are polynomial endomorphisms of $k^n$ then $(F\circ G)^*=G^*\circ F^*$. By definition, $J(F)=J(F^*)$. The chain rule gives that
\bee\label{Yf5}
J(G\circ F)=J(F^*\circ G^*)=F^*(J(G^*))J(F^*)=F^*(J(G))J(F).
\eee

\begin{lm}\label{Yl8} Let $T,S\in \mathscr{L}_n$. Then the following statements are true:

(i) $J(T\cdot S)=T(J(S))+J(S)J(T)$;

(ii) $J([T,S])=T(J(S))-S(J(T))$;

(iii) $\mathrm{div}([T,S])=T(\mathrm{div}(S))-S(\mathrm{div}(T))$.
\end{lm}
\Proof Suppose that $T=D_F$ and $S=D_G$. Then $T\cdot S=D_{D_F(G)}$. Consider the endomorphism $(X+tF)^*$ where $t$ is an independent parameter. Obviously,
\bes
(X+tF)^*(G)=G+tD_F(G)+t^2G_2+\ldots.
\ees
Consequently, $D_F(G)=\frac{\partial}{\partial t} ((X+tF)^*G)|_{t=0}$.
By (\ref{Yf5}), we get
\bes
J((X+tF)^*(G))=J((X+tF)^*\circ G^*) =(X+tF)^*(J(G))J(X+tF)\\
= (J(G)+tD_F(J(G))+t^2T_2+\ldots)(I+tJ(F))\\
= J(G)+t(D_F(J(G))+J(G)J(F))+t^2M_2+\ldots.
\ees
Hence
\bes
J(D_F(G))=\frac{\partial}{\partial t} J((X+tF)^*(G))|_{t=0}=D_F(J(G))+J(G)J(F),
\ees
which proves (i). Notice that (i) directly implies (ii). Besides, $\mathrm{Tr}$ is a linear function and for any
$D\in \mathscr{L}_n$ and $B\in M_n(A)$ we have $\mathrm{Tr}(D(B))=D(\mathrm{Tr}(B))$. Consequently, (ii) implies (iii). $\Box$

\begin{lm}\label{Yl9} Let $D\in\mathscr{L}_n$. Then $J(D)$ is nilpotent if and only if $\mathrm{div}(D^{[q]})=0$ for all $q\geq 1$.
\end{lm}
\Proof By Lemma \ref{Yl8}, we get $J(D^{[2]})=D(J(D))+J(D)^2$ and
\bes
J(D^{[i+1]})=J(D^{[i]}\cdot D)=D^{[i]}(J(D))+J(D)J(D^{[i]}).
\ees
This allows us to prove, by induction on $i$, that
\bee\label{Yf6}
J(D^{[i]})=D^{[i-1]}(J(D))+J(D)D^{[i-2]}(J(D))+\ldots+J(D)^{i-2}D(J(D))+J(D)^i
\eee
for all $i\geq 1$.

Suppose that $J(D)$ is nilpotent. It is well known that $J(D)$ is nilpotent if and only if $\mathrm{Tr}(J(D)^q)=0$ for all $q\geq 1$. Recall that $\mathrm{Tr}(TS)=\mathrm{Tr}(ST)$ for any $T,S\in M_n(A)$. Consequently, for any $D\in \mathscr{L}_n$, $T\in M_n(A)$, and integer $s\geq 1$ we have
\bes
\mathrm{Tr}(D(T^s))
=\mathrm{Tr}(D(T)T^{s-1}+T D(T^2)T^{s-2}
+\\
\ldots+T^{s-2}D(T)T+T^{s-1}D(T))
=s\mathrm{Tr}(T^{s-1}D(T))
\ees
and consequently,
\bee\label{Yf7}
D(\mathrm{Tr}(T^s))=\mathrm{Tr}(D(T^s))
=s\mathrm{Tr}(T^{s-1}D(T))
\eee
Hence $\mathrm{Tr}(T^{s-1}D(T))=0$ and (\ref{Yf6}) implies that $\mathrm{div}(D^{[i]})=\mathrm{Tr}(J(D^{[i]}))=0$.

Suppose that $\mathrm{div}(D^{[q]})=0$ for all $q\geq 1$. We prove by induction on $s$ that $\mathrm{Tr}(J(D)^s)=0$ for all $s\geq 1$. Suppose that it is true for all $s$ such that $1\leq s<i$. Then, (\ref{Yf7}) gives that $\mathrm{Tr}(J(D)^{s-1}D^{[p]}(J(D)))=0$. Consequently, (\ref{Yf7}) implies that $\mathrm{Tr}((J(D)^i)=0$. $\Box$

Let $D$ be an arbitrary element of $\mathscr{L}_n$. Recall that $L_D$ is the Lie algebra generated by all right powers $D^{[i]}$ ($i\geq 1$) of $D$.

\begin{theor}\label{Yt1} Let $D\in\mathscr{L}_n$. Then the Jacobian matrix $J(D)$ of $D$ is nilpotent if and only if the divergence of every element of $L_D$ is zero.
\end{theor}
\Proof This is a direct corollary of Lemmas \ref{Yl8} and \ref{Yl9}. $\Box$

Denote by $I(D)$ the $L_D$-closed subalgebra of $A$ generated by all $\mathrm{Tr}(J(D)^i)=0, i\geq 1$.
\begin{co}\label{Yc4} Let $D\in\mathscr{L}_n$. Then the divergence of every element of $L_D$ belongs to $I(D)$.
\end{co}
\Proof The proof of Lemma \ref{Yl9} can be easily adjusted to prove that $\mathrm{div}(D^{[i]})\in I(D)$. Then Lemma \ref{Yl8} finishes the proof of the corollary. $\Box$

The Lie algebra $L_D$ is a small part of the left-symmetric algebra $\mathscr{L}_D$ generated by $D$. Probably $L_D$ is the maximal subspace of $\mathscr{L}_D$ whose divergence belong to $I(D)$. In other words, I think that if $J(D)$ is nilpotent then $L_D$ is the maximal subspace of elements of $\mathscr{L}_D$ whose divergence are zeroes.

Recall that a derivation $D$ is called {\em  triangular} if $D(x_i)\in k[x_1,\ldots,x_i]$ for all $i$ and {\em strongly triangular} if $D(x_i)\in k[x_1,\ldots,x_{i-1}]$ for all $i$. If $D$ is a triangular derivation with a nilpotent Jacobian matrix $J(D)$, then it is easy to check that $D$ is strongly triangular.
If $D$ is strongly triangular then $J(D)$ is nilpotent and both algebras $L_D$ and $\mathscr{L}_D$ are nilpotent.

{\bf Example 2.}
Now we give an example of derivation $D$ with a nilpotent Jacobian matrix $J(D)$ such that $L_D$ is not nilpotent nor solvable. Consider the automorphism
\bes
(x+s(xt-ys),y+t(xt-ys),s+t^3,t)
\ees
 of the polynomial algebra $k[x,y,s,t]$ studied A. van den Essen \cite{Essen85} and  G. Gorni and G. Zampieri \cite{GZ96}. Put
\bes
F=(s(xt-ys),t(xt-ys),t^3,0).
\ees
Obviously, $J(F)$ is nilpotent. Consider
\bes
D=D_F=s(xt-ys)\partial_x+t(xt-ys)\partial_y+t^3\partial_s.
\ees
 Corollary  \ref{Yc2} gives that $D$ is a right nilpotent element of $\mathscr{L}_n$. Put $w=xt-ys$. Then,
\bes
D(w)=-yt^3, \ \ D(D(w))=-t^4w.
\ees
Consequently, $D$ is not a locally nilpotent derivation and is not a left nilpotent element of  $\mathscr{L}_n$ by Lemma \ref{Yl3}.
Direct calculations give
\bes
D^2=D^{[2]}=t^3(xt-2ys)\partial_x -yt^4\partial_y, \ \
 D^{[2]}(w)=wt^4,\\
  D^{[3]}=st^4w\partial_x+t^5w\partial_y, \ \ D^{[3]}(w)=0, \ \ D^{[4]}=0.
\ees
Consequently, the Lie algebra $L_D$ is generated by two elements $a=D$, $b=D^{[2]}$, and $c=D^{[3]}$. Moreover, we have
\bes
[a,b]= -2c-2A, \ \ A=t^6y\partial_x, \ \  [b,c]=2t^4c, \ \ [a,c]=t^4b.
\ees
These relations show that $L_D$ is not nilpotent. We also have
\bes
[a,A]=t^4b, [A,c]=-t^7b, [A,b]=2t^4 A.
\ees
Let $M$ be the subalgebra of $L_D$ generated by $A,b,c$. Note that $t$ is a constant for all elements of $L_D$.
The homomorphic image of $M$ under $t\mapsto 1$ becomes a Lie algebra with a linear basis $A,b,c$ and and satisfies the relations
\bes
[b,c]=2c,  [A,c]=-b, [A,b]=2 A.
\ees
Consequently, $M$ is not solvable and so is $L_D$.

This example also shows some limits of divergence calculations. The divergence of every element of $L_D$ is zero, but
$L_D$ is not nilpotent nor solvable.

\section{Primitives of the Hopf algebra $\mathrm{NSymm}$}

\hspace*{\parindent}

As an algebra $\mathrm{NSymm}$ \cite{Gelfand} is the free associative algebra
\bes
\mathrm{NSymm}=k\langle Z_1,Z_2,\ldots,Z_n,\ldots\rangle
\ees
 over $k$ in the variables $Z_1,Z_2,\ldots,Z_n,\ldots$. The comultiplication $\bigtriangleup$ and  the counit $\epsilon$ are algebra maps determined by
\bes
\bigtriangleup(Z_n)=\sum_{i+j=n} Z_i\otimes Z_j   \ (Z_0=1), \ \ \ \varepsilon(Z_n)=0,
\ees
for all $n\geq 1$, respectively. The antipod $S$ is an antiisomorphism determined by
\bes
S(Z_n)=\sum_{i_1+\ldots+i_p=n} (-1)^p Z_{i_1}Z_{i_2}\ldots Z_{i_p}
\ees
for all $n\geq 1$.

The Hopf algebra of noncommutative symmetric functions was introduced in \cite{Gelfand} and many systems of free generators and relations between them were described. It was also proved \cite{Gelfand} that
$\mathrm{NSymm}$ is canonically isomorphic to the Solomon descent algebra \cite{Solomon}. It is also known \cite{Gelfand,MR} that the graded dual of $\mathrm{NSymm}$ is the Hopf algebra of quasisymmetric functions $\mathrm{QSymm}$ \cite{Gessel}.

Denote by $\mathrm{Prim}$ the set of all primitive elements  of $\mathrm{NSymm}$, i.e.,
\bes
\mathrm{Prim}= \{p\in \mathrm{NSymm} | \bigtriangleup(p)=p\otimes 1+1\otimes p\}.
\ees
Define the system of elements $U_1,U_2,\ldots,U_i,\ldots$ by
\bes
\sum_{i=1}^{\infty} t^iU_i= \mathrm{log}(\sum_{i=0}^{\infty} t^iZ_i).
\ees
Direct calculations give
\bes
U_m=\sum_{i_1+\ldots+i_k=m} \frac{(-1)^{k-1}}{k} Z_{i_1}\ldots Z_{i_k}
\ees
and
\bes
Z_m=\sum_{i_1+\ldots+i_k=m} \frac{1}{k!} U_{i_1}\ldots U_{i_k}
\ees
for all $m\geq 1$.

It is well known \cite{Gelfand,MR} the Lie algebra $\mathrm{Prim}$ is a free Lie algebra freely generated by $U_1,U_2,\ldots,U_m\ldots$ and $\mathrm{NSymm}$ is the universal enveloping algebra of $\mathrm{NSymm}$.

Consider the following two systems of elements of $\mathrm{NSymm}$ :
\bee\label{Yf8}
\Theta_n(Z)=\sum_{r_1+\ldots+r_k=n} (-1)^{k-1} r_1 Z_{r_1}Z_{r_2}\ldots Z_{r_k},
\eee
and
\bee\label{Yf9}
\Psi_n(Z)=\sum_{r_1+\ldots+r_k=n} (-1)^{k-1} r_k Z_{r_1}Z_{r_2}\ldots Z_{r_k},
\eee
where $r_i\in \mathbb{N}=\{1,2,\ldots\}$ and $n\geq 1$.

Notice that in our notations, $Z_i$ correspond to complete symmetric functions $S_i$, $\Psi_i$ are the power sums symmetric functions, and $U_i$ correspond to power sums of the second kind $\Phi_i/i$ in \cite{Gelfand}. The functions corresponding to $\Theta_i$ were not considered in \cite{Gelfand} since $\Theta_i$ can be obtained from $\Psi_i$ by the natural involution of $\mathrm{NSymm}$ preserving all $Z_i$. But in needs of the Jacobian Conjecture it is necessary to study the relations between $\Theta_i$ and $\Psi_i$ more deeply.

The systems of elements (\ref{Yf8}) and (\ref{Yf9}) are primitive systems of free generators of the free associative algebra $\mathrm{NSymm}$ \cite{Gelfand} and can be defined recursively by
\bes
nZ_n=\Theta_n(Z)+\Theta_{n-1}Z_1+\Theta_{n-2}Z_2+\ldots+\Theta_1Z_{n-1}
\ees
and
\bes
nZ_n=\Psi_n(Z)+Z_1 \Psi_{n-1}+Z_2 \Psi_{n-2}+\ldots+Z_{n-1}\Psi_1
\ees
for all $n\geq 1$.

Recall that a composition is a vector $I = (i_1,\ldots,i_m)$ of nonnegative integers, called
the parts of $I$. The length $l(I)$ of the composition $I$ is the number $k$ of its parts and the
weigt of $I$ is the sum $|I| =\Sigma i_j$ of its parts. We use notations
\bes
Z^I=Z_{i_1}\ldots Z_{i_m}, \ \ \  \Theta^I=\Theta_{i_1}\ldots \Theta_{i_m}, \ \ \  \Psi^I=\Psi_{i_1}\ldots \Psi _{i_m}.
\ees
Put also
\bes
\pi_u(I)=i_1(i_1+i_2)\ldots (i_1+i_2+\ldots+i_m)
\ees
and $\mathrm{lp}(I)=i_m$ (the last part of $I$). Let $J$ be another composition. We say that $I\preceq J$ if
$J=(J_1,\ldots,J_m)$ and $|J_j|=i_j$ for all $j$. For example, $(3,2,6)\preceq (2,1,2,3,1,2)$. If $I\preceq J$ then put
\bes
\pi_u(J,I)=\prod_{i=1}^m \pi_u(J_i),\ \ \  \mathrm{lp}(J,I)=\prod_{i=1}^m \mathrm{lp}(J_i).
\ees
The following formulas are proved in \cite{Gelfand}.
\bee\label{Yf10}
Z^I=\sum_{J\succeq I} \ \frac{1}{\pi_u(J,I)} \Psi^J, \ \ \  \Psi^I=\sum_{J\succeq I} (-1)^{l(J)-l(I)} \ \mathrm{lp}(J,I) Z^J.
\eee
Denote by $w$ the natural involution of the free associative algebra $\mathrm{NSymm}$ preserving all $Z_i$. Obviously, $w(\Theta_i)=\Psi_i$ and $w(\Psi_i)=\Theta_i$ for all $i$. Applying $w$, from (\ref{Yf10}) we get
\bes
Z^{\overline{I}}=\sum_{J\succeq I} \ \frac{1}{\pi_u(J,I)} \Theta^{\overline{J}}, \ \ \  \Theta^{\overline{I}}=\sum_{J\succeq I} (-1)^{l(J)-l(I)} \ \mathrm{lp}(J,I) Z^{\overline{I}},
\ees
where $\overline{I}$ is the mirror image of the composition $I$, i.e. the new composition
obtained by reading $I$ from right to left.

Consequently,
\bee\label{Yf11}
Z^I=\sum_{J\succeq I} \ \frac{1}{\pi_u(\overline{J},\overline{I})} \Theta^J, \ \ \  \Theta^I=\sum_{J\succeq I} (-1)^{l(J)-l(I)} \ \mathrm{lp}(\overline{J},\overline{I}) Z^J.
\eee
Using (\ref{Yf9}) and (\ref{Yf11}), we get
\bes
\Psi_n=\sum_{|I|=n} (-1)^{l(I)-1} \ \mathrm{lp}(I) Z^I=\sum_{|I|=n} (-1)^{l(I)-1} \ \mathrm{lp}(I) \sum_{J\succeq I} \ \frac{1}{\pi_u(\overline{J},\overline{I})} \Theta^J,
\ees
i.e.,
\bee\label{Yf12}
\Psi_n=\sum_{J\succeq I, |I|=n} (-1)^{l(I)-1}   \frac{\mathrm{lp}(I)}{\pi_u(\overline{J},\overline{I})} \Theta^J.
\eee
In fact, $\Psi_n$ can be expressed as a Lie polynomial of $\Theta_1,\ldots,\Theta_n$. We have
\bes
\Psi_1=\Theta_1, \ \Psi_2=\Theta_2, \ \Psi_3= \Theta_3+1/2 [\Theta_2,\Theta_1],\\
 \Psi_4=\Theta_4+2/3[\Theta_3,\Theta_1]+1/6[[\Theta_2,\Theta_1],\Theta_1].
\ees
It will be interesting to find the Lie expression of $\Psi_n$ in $\Theta_1,\ldots,\Theta_n$.

\section{An action of the Hopf algebra $\mathrm{NSymm}$}

\hspace*{\parindent}

We define an action
\bes
\mathrm{NSymm} \times P_n \longrightarrow P_n  \ \ \ \ ((T,a)\mapsto T\circ a)
\ees
of $\mathrm{NSymm}$ on the polynomial algebra $P_n$ related to an $n$-tuple $F$. Since  $\mathrm{NSymm}$ is a free associative algebra, it is sufficient to define $Z_i\circ a$ for all $i\geq 1$ and $a\in P_n$. For any $a\in P_n$ there exists a unique system of elements $g_i\in P_n, i\geq 1$ such that
\bes
(X+tF)^*(a)=a+tg_1+t^2g_2+\ldots+t^ng_n+\ldots,
\ees
where $t$ is an independent variable. Put $Z_i\circ a=g_i$ for all $i\geq 1$. Then
\bes
(X+tF)^*(a)=a+tZ_1(a)+t^2Z_2(a)+\ldots+t^nZ_n(a)+\ldots.
\ees
This formula can be considered as a linearization of the action of $(X+tF)^*$ on $P_n$.
 Denote by
\bes
\lambda : \mathrm{NSymm} \longrightarrow \mathrm{Hom}_k(P_n,P_n)
\ees
the homomorphism corresponding to this representation, where $\mathrm{Hom}_k(P_n,P_n)$ is the set of all $k$-linear maps from $P_n$ to $P_n$. First of all we show that $\lambda(\mathrm{NSymm})\subseteq A_n$.

Denote by $p : P_n\otimes_k P_n\rightarrow P_n$ the product in the polynomial algebra $P_n$.
\begin{lm}\label{Yl10}  Let $T\in \mathrm{NSymm}$. Then $\lambda(T)p=p\lambda(\bigtriangleup(T))$.
\end{lm}
\Proof It is easy to check that the set of elements $T\in \mathrm{NSymm}$ satisfying the statement of the lemma forms a subalgebra. Consequently, we may assume that $T=Z_n$. If $a,b\in A$ then
\bes
\sum_{i=0} t^i\lambda(Z_i)(ab)=(X+tF)^*(ab)\\
=((X+tF)^*(a))((X+tF)^*(b))\\
=(\sum_{i=0} t^i\lambda(Z_i)(a)) (\sum_{i=0} t^i\lambda(Z_i)(b)).
\ees
Comparing coefficients in the degrees of $t$ we get $Z_i(ab)=\sum_{r+s=i} Z_r(a)Z_s(b)$. This means $Z_ip=p\bigtriangleup(Z_i)$. $\Box$
\begin{lm}\label{Yl11}  $\lambda(\mathrm{Prim})\subseteq W_n$ and $\lambda(\mathrm{NSymm})\subseteq A_n$.
\end{lm}
\Proof If $T\in \mathrm{Prim}$ then, by Lemma \ref{Yl10}, we get
\bes
\lambda(T)(ab)=\lambda(T)p(a\otimes b)=p\lambda(\bigtriangleup(T))(a\otimes b)
= p\lambda(T\otimes 1+ 1\otimes T)(a\otimes b)\\
=p(\lambda(T)(a)\otimes b+a\otimes \lambda(T)(b))
=\lambda(T)(a)b+a\lambda(T)(b),
\ees
i.e., $\lambda(T)\in W_n$.

Notice that $\mathrm{NSymm}$ is a free associative algebra and any action of $\mathrm{NSymm}$ is well defined by the action of any free system of generators. For example, $\Theta_1,\Theta_2,\ldots,\Theta_n,\ldots\in \mathrm{Prim}$ is a free system of generators of $\mathrm{NSymm}$ and $\lambda(\Theta_1),\lambda(\Theta_2),\ldots,\lambda(\Theta_n),\ldots\in W_n$. Consequently, for any $T\in \mathrm{NSymm}$ element $\lambda(T)$ is a differential operator on $P_n$, i.e., $\lambda(T)\in A_n$.
$\Box$

By this lemma, we have a homomorphism
\bee\label{Yf13}
\lambda : \mathrm{NSymm} \longrightarrow A_n.
\eee

\begin{lm}\label{Yl12}
Let $a\in P_n$ and $\deg\,a\leq k$. Then $\lambda(Z_i)(a)=0$ for all $i\geq k+1$.
\end{lm}
\Proof Obviously, the degree of $(X+tF)(a)=a((x_1+tf_1),\ldots,(x_n+tf_n))$ with respect to $t$ is less than or equal to $k$. Consequently, $\lambda(Z_i)(a)=0$ for all $i\geq k+1$. $\Box$

\begin{pr}\label{Ypr1} Let
\bes
(X+tF)^{-1}=X+tF_1+t^2F_2+\ldots+t^mF_m+\ldots
\ees
be the formal inverse to the endomorphism $X+tF$ of $k[t]^n$. Then
\bes
-\lambda(\Psi_m)(X)=F_m
\ees
for all $m\geq 1$.
\end{pr}
\Proof Consider the endomorphism $(X+tF)^* : k[t]\otimes_k P_n\rightarrow k[t]\otimes_k P_n$ of the $k[t]$-algebra.  Notice that
\bee
(X+tF)^*=1+t\lambda(Z_1)+t^2\lambda(Z_2)+\ldots+t^n\lambda(Z_n)+\ldots
\eee
by the definition of $\lambda(Z_i)$. Then,
\bes
(X+tF)^*=1-T, \ \ T=-(t\lambda(Z_1)+t^2\lambda(Z_2)+\ldots+t^n\lambda(Z_n)+\ldots),
\ees
and
\bes
((X+tF)^*)^{-1}=1+T+T^2+\ldots +T^n+\ldots.
\ees
Direct calculation gives
\bes
((X+tF)^*)^{-1}=1+tT_1+t^2T_2+\ldots+t^nT_n+\ldots,
\ees
where
\bes
T_m=\sum_{r_1+\ldots+r_k=m} (-1)^k\lambda(Z_{r_1})\lambda(Z_{r_2})\ldots \lambda(Z_{r_k}), \ \ n\geq 1.
\ees
Notice that $(X+tF)^{-1}=((X+tF)^*)^{-1}(X)$ and $F_m=T_m(X)$. Then,
\bes
F_n=\sum_{r_1+\ldots+r_k=m} (-1)^k\lambda(Z_{r_1})\lambda(Z_{r_2})\ldots \lambda(Z_{r_k})(X)\\
=\sum_{r_1+\ldots+r_k=m} (-1)^kr_k \lambda(Z_{r_1})\lambda(Z_{r_2})\ldots \lambda(Z_{r_k})(X)
\ees
by lemma \ref{Yl12}. Consequently, $F_m=-\lambda(\Psi_m)(X)$. $\Box$

The homomorphism (\ref{Yf13}) coincides with one of a series of homomorphisms constructed in \cite{Zhao08} and the images of primitive generators were calculated in \cite{Zhao08}.
\begin{lm}\label{Yl13} $(i)$ \ $\lambda(\Theta_m)=(-1)^{m-1}D_F^{[m]}$ for all $m\geq 1$.

$(ii)$ \ $\lambda(\Psi_m)=-D_{F_m}$ for all $m\geq 1$.
\end{lm}
\Proof By Lemma \ref{Yl11}, $\lambda(\Theta_m)$ and $\lambda(\Psi_m)$ are derivations of $P_n$. Consequently, it is sufficient to prove that $\lambda(\Theta_m)(X)=(-1)^{m-1}D_F^{[m]}(X)$ and $\lambda(\Psi_m)(X)=-D_{F_m}(X)=-F_m$. Proposition \ref{Ypr1} implies $(ii)$. We have $\lambda(\Theta_1)(X)=F=D_F(X)$ since $\Theta_1=Z_1$. Then, $\lambda(\Theta_1)=D_F$. Leading an induction on $m$, by (\ref{Yf8}) and Lemma \ref{Yl12}, we get
\bes
\lambda(\Theta_m)(X)=-\lambda(\Theta_{m-1})\lambda(Z_1)(X)=-\lambda(\Theta_{m-1})\lambda(Z_1)(X)\\
=(-1)^{m-1}D_F^{[m-1]}(F)=(-1)^{m-1}D_F^{[m]}(X). \ \ \ \Box
\ees

Put $D=D_F$. Recall that $L_D$ is the subalgebra of $W_n$ generated by all right powers $D^{[m]}$ ($m\geq 1$) of $D$ and $A_D$ is the subalgebra of $A_n$ generated by the same elements.
\begin{co}\label{Yc5} Let $D=D_F$. Then $\lambda(\mathrm{Prim})=L_D$ and $\lambda(\mathrm{NSymm})=A_D$.
\end{co}
\Proof This is an immediate corollary of Lemmas \ref{Yl11} and \ref{Yl13}.
$\Box$

\begin{theor}\label{Yt2} Let $F=(f_1,\ldots,f_n)$ be an arbitrary $n$-tuple of the polynomial algebra $P_n=k[x_1,\ldots,x_n]$, $L_D=L_F$ be the Lie algebra generated by all right powers $D^{[m]}$ ($m\geq 1$) of $D=D_F$, and
\bes
(X+tF)^{-1}=X+tF_1+t^2F_2+\ldots+t^mF_m+\ldots
\ees
be the formal inverse to the endomorphism $X+tF$ of $k[t]^n$. Then the Lie algebra $L_D$ is generated by all $D_{F_m}$ where $m\geq 1$.
\end{theor}
\Proof By Corollary \ref{Yc5}, $L_D$ is the image of the Lie algebra $\mathrm{Prim}$ of all primitive elements of $\mathrm{NSymm}$ under $\lambda$. The set of elements $\Psi_m$, where $m\geq 1$, is also generates $\mathrm{Prim}$.  Consequently, Lemma \ref{Yl13} implies the statement $(i)$. $\Box$

One more interesting system of generators $\lambda(U_1),\ldots,\lambda(U_m),\ldots$ of the Lie algebra $L_D$ corresponds to the coefficients of $D-\mathrm{log}$ of $X+tF$ considered in \cite{WZh,Zhao08}.

\section{Comments and some open questions}

\hspace*{\parindent}

So, we introduced three algebras $A_F$, $L_F$, and $\mathscr{L}_F$ related to the study of the Jacobian Conjecture, i.e., to the study of the polynomial endomorphism $X+tF$ with a nilpotent Jacobian matrix $J(F)$. If $J(F)$ is nilpotent then $D=D_F$ is right nilpotent by Corollary \ref{Yc2}. Let $p$ be a positive integer such that $D^{[p]}=0$. In order to solve the Jacobian Conjecture, it is necessary to prove that there exists $m=m(F)$ such that $F_i=0$ for all $i\geq m$ in notations of Theorem \ref{Yt2}.  Using Lemma \ref{Yl13} and (\ref{Yf12}), we get
\bee\label{Yf15}
D_{F_i}=\sum_{p\geq J\succeq I, |I|=n} (-1)^{l(I)+|J|-l(J)}   \frac{\mathrm{lp}(I)}{\pi_u(\overline{J},\overline{I})} D^J,
\eee
where $D^J=D^{[j_1]}\ldots D^{[j_s]}$ for any $J=(j_1,\ldots,j_s)$ and $p\geq J$ means that $p\geq j_i$ for all $i$. Moreover, the right hand side of this equation is a Lie polynomial in $D^{[s]}$ where $s\geq 1$ and the Jacobian Conjecture can be considered as a problem of the algebra $L_F$. But I cannot see how to use the degree of $F$ in this formula. We cannot prove that $D_{F_i}=0$ without this.

Let's come back to formula (\ref{Yf10}) and Lemma \ref{Yl12}. Suppose that the degree of $F$ is $m$. A composition
$I=(i_k,\ldots,i_1)$ of length $k$ is called $m$-reduced if $i_1=1$, $i_2\leq m$, and $i_j\leq (i_1+\ldots+i_{j-1})(m-1)+1$ for all $3\leq j\leq k$. Let $T_n$ be the set of all $m$-reduced compositions $I$ with $|I|=n$. Notice that $\mathrm{lp}(I)=1$ if $I$ is $m$-reduced. If $I$ is not $m$-reduced then $\lambda(Z^I)(X)=0$ by Lemma \ref{Yl12}. For this reason we can consider only $m$-reduced compositions in (\ref{Yf10}). Then we get
\bee\label{Yf16}
D_{F_i}=\sum_{J\succeq I, I\in T_n} (-1)^{l(I)+|J|-l(J)}   \frac{1}{\pi_u(\overline{J},\overline{I})} D^J
\eee
in $A_F$ but not in $L_F$. So, we did not get $D_{F_i}=0$ yet. In fact, to derive (\ref{Yf16}) we used only the nilpotency of $D$ and the degree of $F$. In connection with this, the following question is very interesting.
\begin{prob}\label{Yp1} Is the Jacobian matrix $J(D)$ of $D$ nilpotent if $D$ is a right nilpotent element of $\mathscr{L}_n$?
\end{prob}
If the answer to this question is negative, then we probably cannot prove that $D_{F_i}=0$ in $A_F$.

The formula (\ref{Yf16}) can be considered as a formula in the left-symmetric algebra $\mathscr{L}_D$ where the associative product $D^J$ is changed by the left normed product. For this reason  left operator identities of $\mathscr{L}_n$ are very important. Notice that $J(F)$ is nilpotent if and only if $R_D$ is nilpotent by Corollary \ref{Yc2}. So, this condition is expressed in the language of right multiplication operators but (\ref{Yf16}) is expressed in the language of left operators.

The following problem is interesting in connection with Lemmas \ref{Yl4} and \ref{Yl6}.
\begin{prob}\label{Yp2} Describe the structure of the multiplication algebra $M(\mathscr{L}_n)$ of the left-symmetric algebra $\mathscr{L}_n$.
\end{prob}
It is well known that all trace identities of matrix algebras are corollaries of the Cayley-Hamilton trace identities \cite{Razmyslov74}.
\begin{prob}\label{Yp3} Is every trace identity (or identity) of $\mathscr{L}_n$ a corollary of the Cayley-Hamilton trace identities (\ref{Yf4}).
\end{prob}
By Lemma \ref{Yl7}, a positive answer to this question implies that every identity of $W_n$ is a corollary of the
Cayley-Hamilton trace identities.

In order to solve the Jacobian Conjecture we need more information about left operator identities of $\mathscr{L}_n$.
\begin{prob}\label{Yp4} Describe all left operator identities of $\mathscr{L}_n$.
\end{prob}

It is interesting to know that what types of properties can be better described in the language of $A_D$.
\begin{prob}\label{Yp5} Describe all $D\in \mathscr{L}_n$ such that $A_D$ is a simple algebra.
\end{prob}

\begin{prob}\label{Yp6} Is there any derivation $D$ with nilpotent Jacobian matrix $J(D)$ such that $A_D$ is a simple algebra?
\end{prob}

Example 1 shows that the nilpotency of $J(F)$ does not imply neither nilpotency nor solvability of $L_F$.
\begin{prob}\label{Yp7} Describe necessary and sufficient conditions of the nilpotency (and solvability) of the Lie algebra $L_D$.
\end{prob}

 At the moment I know that $\mathscr{L}_D$ is nilpotent if and only if $\mathrm{div}(\mathscr{L}_D)=0$.

\bigskip

\begin{center}
{\bf\large Acknowledgments}
\end{center}

\hspace*{\parindent}

I am grateful to Max-Planck Institute f\"ur Mathematik for their
hospitality and excellent working conditions, where part of this work has been done.

\end{document}